# Laws of the single logarithm for delayed sums of random fields

ALLAN GUT[1]   and ULRICH STADTMÜLLER[2]

[1]*Department of Mathematics, Uppsala University, Box 480, SE-751 06 Uppsala, Sweden.
E-mail: allan.gut@math.uu.se, url: http://www.math.uu.se/˜allan/*
[2]*Department of Number Theory and Probability Theory, Ulm University, D-89069 Ulm,
Germany. E-mail: ulrich.stadtmueller@uni-ulm.de,
url: http://www.mathematik.uni-ulm.de/matheIII/members/stadtmueller/stadtmueller.html*

We extend a law of the single logarithm for delayed sums by Lai to delayed sums of random fields.
A law for subsequences, which also includes the one-dimensional case, is obtained in passing.

*Keywords:* delayed sums; law of the iterated logarithm; law of the single logarithm;
multidimensional indices; random fields; sums of i.i.d. random variables; window

## 1. Introduction

The aim of this paper is to study certain strong laws of large numbers for random fields
in the sense that we shall discuss the (possible) equivalence of certain limit relations for
sums over random fields of i.i.d. random variables with suitable moment conditions.

More precisely, let $\{X_{\mathbf{k}}, \mathbf{k} \in \mathbf{Z}_+^d\}$ be i.i.d. random variables with partial sums $S_{\mathbf{n}} = \sum_{\mathbf{k} \le \mathbf{n}} X_{\mathbf{k}}$, $\mathbf{n} \in \mathbf{Z}_+^d$, where the random field or index set $\mathbf{Z}_+^d$, $d \ge 2$, is the positive integer
$d$-dimensional lattice with coordinatewise partial ordering $\le$.

There exist various results on strong laws for random fields in the literature. We now
very briefly describe some of them.

For random fields with i.i.d. random variables $\{X_{\mathbf{k}}, \mathbf{k} \in \mathbf{Z}_+^d\}$, the analog of Kol-
mogorov's strong law (see [17]) reads as follows:

$$\frac{S_{\mathbf{n}}}{|\mathbf{n}|} = \frac{1}{|\mathbf{n}|} \sum_{\mathbf{k} \le \mathbf{n}} X_{\mathbf{k}} \xrightarrow{\text{a.s.}} 0 \quad \Longleftrightarrow \quad \mathrm{E}(|X|(\log^+ |X|)^{d-1}) < \infty, \qquad \mathrm{E}X = 0.$$

Here and throughout, $\log^+ x = \max\{\log x, 1\}$.

In [7], the first author has considered Marcinkiewicz laws for i.i.d. random fields. To be more precise he proved that, for $0 < r < 2$,

$$\frac{1}{|\mathbf{n}|^{1/r}} S_{\mathbf{n}} \overset{\text{a.s.}}{\to} 0 \quad \Longleftrightarrow \quad \mathrm{E}(|X|^r (\log^+ |X|)^{d-1}) < \infty \qquad (\text{and } \mathrm{E}X = 0 \text{ when } r \geq 1).$$

Under somewhat stronger moment conditions, strong laws can be proven under relaxed conditions on the distribution and the dependence structure. It has, for example, been shown that for $d$-dimensional martingales* (see, e.g., [14, 19] for a definition),

$$\frac{S_{\mathbf{n}}}{|\mathbf{n}|} \overset{\text{a.s.}}{\to} 0 \qquad \text{if } \sum_{\mathbf{n}} \frac{\mathrm{E}|S_{\mathbf{n}}|^{2q}}{|\mathbf{n}|^{q+1}} < \infty \qquad \text{for some } q \geq 1.$$

For a proof, see [14], where orthogonal random fields were also discussed.

As for the rate of convergence in the strong law for random fields for i.i.d. random variables, there exists, as in the one-dimensional case, a law of the iterated logarithm (LIL), which reads as follows in the case $d > 1$ (see [22]):

$$\limsup_{\mathbf{n} \to \infty} \left( \liminf_{\mathbf{n} \to \infty} \right) \frac{S_{\mathbf{n}}}{\sqrt{2|\mathbf{n}| \log \log |\mathbf{n}|}} = \sigma \sqrt{d} (-\sigma \sqrt{d}) \qquad \text{a.s.} \quad \Longleftrightarrow$$

$$\mathrm{E}\left( X^2 \frac{(\log^+ |X|)^{d-1}}{\log^+ \log^+ |X|} \right) < \infty \quad \text{and} \quad \mathrm{E}X = 0, \qquad \mathrm{E}X^2 = \sigma^2.$$

In particular, the moment condition and the limit depend on the dimension $d$. If one restricts the $\limsup$ to a sector avoiding the boundaries of $\mathbf{Z}_+^d$, the first author has proven [9] that the law of the iterated logarithm then holds under the same moment condition as in the case $d = 1$ and has limit points $\pm \sigma$.

A different question concerns the limit behaviour for delayed sums, sometimes called *lag sums*. Let us start with ordinary random variables $\{X_k, k \geq 1\}$ and let $T_{n,n+n^\alpha} = S_{n+n^\alpha} - S_n = \sum_{k=n+1}^{n+n^\alpha} X_k$, where $0 < \alpha < 1$ and, to be precise, $n^\alpha := [n^\alpha]$. For i.i.d. summands $\{X_k, k \geq 1\}$, we then have (see [5, 15])

$$\lim_{n \to \infty} \frac{T_{n,n+n^\alpha}}{n^\alpha} = 0 \qquad \text{a.s.} \quad \Longleftrightarrow \quad \mathrm{E}|X|^{1/\alpha} < \infty, \qquad \mathrm{E}X = 0 \quad \text{and}$$

$$\limsup_{n \to \infty} \frac{T_{n,n+n^\alpha}}{\sqrt{2n^\alpha \log n}} = \sqrt{1-\alpha} \qquad \text{a.s.} \quad \Longleftrightarrow \quad \mathrm{E}(|X|^{2/\alpha} (\log^+ |X|)^{-1/\alpha}) < \infty,$$

$$\mathrm{E}X^2 = 1, \qquad \mathrm{E}X = 0.$$

We call the latter result a *law of the single logarithm* (LSL). Such results are of particular interest since they help to evaluate weighted sums of i.i.d. random variables for certain classes of weights. Especially in the case $\alpha = 1/2$, it has been shown (see, e.g., [1, 5]) that certain summability methods, such as those defined by the delayed sums, Euler, Borel and certain Valiron methods, are equivalent for sequences $\{s_n, n \geq 1\}$ satisfying $s_n = o(n^{1/2})$.

The latter remains a.s. true for random variables with the moment condition $EX^2 < \infty$ since the partial sums then satisfy $S_n = o_{a.s.}(n^{1/2})$. Other equivalences with limit relations for delayed sums were given in [3]. These results allow one to prove a law of large numbers and, using a similar idea, an LSL for the associated weighted sums (see, e.g., [4]). Since multivariate summability methods have attracted new interest, multiindex versions of Chow's and Lai's results are of interest.

Hence, the aim of the present paper is to investigate the LSL problem for delayed sums of random fields with i.i.d. components. The strong law itself was established in [20].

The paper is organized as follows. We begin, in Section 2, with definitions and the statement of our main result, after which we collect various preliminaries in Section 3 and the truncation procedure and the Kolmogorov exponential bounds in Section 4. Thus prepared, we present the proof of our main result, the LSL for delayed sums. Examining the proof, it turns out that one can also prove an LSL for subsequences, the result of which is given in Section 6; in particular, this result is also valid for the case $d = 1$. A closing section contains some additional results and remarks.

## 2. Setting and main result

Let $\mathbf{Z}_+^d$, $d \geq 2$, denote the positive integer $d$-dimensional lattice with coordinatewise partial ordering $\leq$, that is, for $\mathbf{m} = (m_1, m_2, \ldots, m_d)$ and $\mathbf{n} = (n_1, n_2, \ldots, n_d)$, $\mathbf{m} \leq \mathbf{n}$ means that $m_k \leq n_k$ for $k = 1, 2, \ldots, d$. Similarly, $\mathbf{n}^\alpha = (n_1^\alpha, n_2^\alpha, \ldots, n_d^\alpha)$. The "size" of a point equals $|\mathbf{n}| = \prod_{k=1}^d n_k$. Moreover, $\mathbf{n} \to \infty$ means that $n_k \to \infty$ for all $k = 1, 2, \ldots, d$. We shall also abuse notation for simplicity and treat the coordinates of $\mathbf{n}^\alpha$ as integers. Finally, $C$ denotes a numerical constant which may change from appearance to appearance.

Throughout the paper, $X$ and $\{X_\mathbf{k}, \mathbf{k} \in \mathbf{Z}_+^d\}$ are i.i.d. random variables with partial sums $S_\mathbf{n} = \sum_{\mathbf{k} \leq \mathbf{n}} X_\mathbf{k}$, $\mathbf{n} \in \mathbf{Z}_+^d$.

For $d = 1$, (the forward) *delayed sums* are $T_{n,n+k} = \sum_{j=n+1}^{n+k} X_j, k \geq 1$, that is, the increments from $S_n$ to $S_{n+k}$. For $d = 2$, the analog is the incremental rectangle

$$T_{\mathbf{n},\mathbf{n}+\mathbf{k}} = S_{n_1+k_1,n_2+k_2} - S_{n_1+k_1,n_2} - S_{n_1,n_2+k_2} + S_{n_1,n_2}$$

and for higher dimensions, the analogous $d$-dimensional cube.

The aim of this paper is to prove a law of the single logarithm for the family of delayed increments or windows

$$\{T_{\mathbf{n},\mathbf{n}+\mathbf{n}^\alpha}, \mathbf{n} \in \mathbf{Z}_+^d\}.$$

Our main result is the following.

**Theorem 2.1.** *Suppose that $\{X_\mathbf{k}, \mathbf{k} \in \mathbf{Z}_+^d\}$ are i.i.d. random variables with mean 0 and finite variance $\sigma^2$, and set $S_\mathbf{n} = \sum_{\mathbf{k} \leq \mathbf{n}} X_\mathbf{k}$, $\mathbf{n} \in \mathbf{Z}_+^d$. If*

$$EX^{2/\alpha}(\log^+ |X|)^{d-1-1/\alpha} < \infty, \tag{2.1}$$

*where $0 < \alpha < 1$, then*

$$\limsup_{\mathbf{n}\to\infty}\left(\liminf_{\mathbf{n}\to\infty}\right)\frac{T_{\mathbf{n},\mathbf{n}+\mathbf{n}^\alpha}}{\sqrt{2|\mathbf{n}|^\alpha\log|\mathbf{n}|}} = +\sigma\sqrt{1-\alpha}\ \ (-\sigma\sqrt{1-\alpha})\qquad a.s. \tag{2.2}$$

*Conversely, if*

$$P\left(\limsup_{\mathbf{n}\to\infty}\frac{|T_{\mathbf{n},\mathbf{n}+\mathbf{n}^\alpha}|}{\sqrt{|\mathbf{n}|^\alpha\log|\mathbf{n}|}} < \infty\right) > 0, \tag{2.3}$$

*then* (2.1) *holds,* $\mathrm{E}X = 0$ *and* (2.2) *holds with* $\sigma^2 = \mathrm{Var}\,X$.

## 3. Preliminaries

We first observe that a partial sum $S_{\mathbf{n}}$ is a sum of $|\mathbf{n}|$ i.i.d. random variables, which implies that distributional properties and various inequalities do not depend on the (partial) order of the index set and thus remain valid "automatically". However, the Lévy inequalities, for example, concern the distribution of $\max_{\mathbf{k}\leq\mathbf{n}} S_{\mathbf{k}}$ and here the structure of the index set enters.

**Lemma 3.1.** *Suppose that* $\{X_{\mathbf{k}}, \mathbf{k} \in \mathbf{Z}_+^d\}$ *are independent random variables with mean 0 and partial sums* $S_{\mathbf{n}} = \sum_{\mathbf{k}\leq\mathbf{n}} X_{\mathbf{k}}$.

(a) *If, in addition, the summands are symmetric, then*

$$P\left(\max_{\mathbf{k}\leq\mathbf{n}} S_{\mathbf{k}} > x\right) \leq 2^d P(S_{\mathbf{n}} > x).$$

(b) *If the variances are finite, then*

$$P\left(\max_{\mathbf{k}\leq\mathbf{n}} S_{\mathbf{k}} > x\right) \leq 2^d P(S_{\mathbf{n}} > x - d\sqrt{2\,\mathrm{Var}(S_{\mathbf{n}})}).$$

The proof is based on induction over the dimensions. For details (in the i.i.d. case), see [8], Lemma 2.3. Two-sided versions are immediate.

We also need relations between tail probabilities and moments analogous to the one-dimensional

$$\mathrm{E}|X| < \infty \quad\iff\quad \sum_{n=1}^\infty P(|X| > n) < \infty.$$

More precisely, we wish to find the necessary moment condition to ensure that

$$\sum_{\mathbf{n}} P(|X| > |\mathbf{n}|) < \infty. \tag{3.1}$$

For this, it turns out that the quantities

$$d(j) = \mathrm{Card}\{\mathbf{k} : |\mathbf{k}| = j\} \quad \text{and} \quad M(j) = \mathrm{Card}\{\mathbf{k} : |\mathbf{k}| \leq j\}$$

and their asymptotics

$$\frac{M(j)}{j(\log j)^{d-1}} \to \frac{1}{(d-1)!} \qquad \text{as } j \to \infty \tag{3.2}$$

and

$$d(j) = o(j^\delta) \qquad \text{for any } \delta > 0 \text{ as } j \to \infty \tag{3.3}$$

play a crucial role. We refer [12], Chapter XVIII and [21], relation (12.1.1) (for the case $d = 2$). The quantity $d(j)$ itself has no pleasant asymptotics in the sense that $\liminf_{j \to \infty} d(j) = d$ and $\limsup_{j \to \infty} d(j) = +\infty$.

We shall also exploit the fact that all terms in expressions such as the sum in (3.1) with equisized indices are equal, which implies that

$$\sum_{\mathbf{n}} P(|X| > |\mathbf{n}|) = \sum_{j=1}^{\infty} \sum_{|\mathbf{n}|=j} d(j) P(|X| > j) \tag{3.4}$$

in that particular case.

This fact, partial summation and (3.2) yield the first part of the following lemma; see also [17]. The second part is a consequence of the fact that the inverse of the function $y = x^\alpha (\log x)^\kappa$ behaves asymptotically like $x = y^{1/\alpha} (\log y)^{-(\kappa/\alpha)}$ (except for some constant factor(s)).

**Lemma 3.2.** *Let $\alpha > 0$ and $\kappa \in \mathbb{R}$ and suppose that $\{X_\mathbf{k}, \mathbf{k} \in \mathbf{Z}_+^d\}$ are i.i.d. random variables with mean 0 and partial sums $S_\mathbf{n} = \sum_{\mathbf{k} \leq \mathbf{n}} X_\mathbf{k}$. Then*

$$\sum_{\mathbf{n}} P(|X| > |\mathbf{n}|^\alpha (\log |\mathbf{n}|)^\kappa) < \infty \quad \Longleftrightarrow \quad \mathrm{E}|X|^{1/\alpha} (\log^+ |X|)^{d-1-\kappa/\alpha} < \infty.$$

For purely numerical sequences, we have the following

**Lemma 3.3.** *Let $\kappa \geq 1$, $\theta > 0$ and $\eta \in \mathbb{R}$.*

$$\sum_{i=2}^{\infty} \sum_{\{\mathbf{n} : |\mathbf{n}| = i^\kappa (\log i)^\eta\}} \frac{1}{|\mathbf{n}|^\theta} = \sum_{i=2}^{\infty} \frac{d(i^\kappa (\log i)^\eta)}{i^{\kappa\theta} (\log i)^{\eta\theta}} \quad \begin{cases} < \infty, & \text{when } \theta > \dfrac{1}{\kappa}, \\[2mm] = \infty, & \text{when } \theta < \dfrac{1}{\kappa}. \end{cases}$$

**Proof.** Recalling (3.4), the convergence part follows from (3.3) and the divergence part from the fact that $d(i) \geq d$ for all $i$. $\qquad\square$

# 4. Truncation and exponential bounds

The typical pattern in proving results of the LIL-type requires two truncations; the first to match the Kolmogorov exponential bounds (see, e.g., [11], Section 8.2) and the second to match the moment requirements.

To this end, let $\delta$ be small, let

$$b_{\mathbf{n}} = b_{|\mathbf{n}|} = \frac{\sigma\delta}{\varepsilon} \frac{\sqrt{|\mathbf{n}|^{\alpha}}}{\log |\mathbf{n}|} \qquad (4.1)$$

and set

$$X'_{\mathbf{n}} = X_{\mathbf{n}} I\{|X_{\mathbf{n}}| \leq b_{\mathbf{n}}\}, \qquad X''_{\mathbf{n}} = X_{\mathbf{n}} I\{b_{\mathbf{n}} < |X_{\mathbf{n}}| < \delta\sqrt{|\mathbf{n}|^{\alpha} \log |\mathbf{n}|}\},$$

$$X'''_{\mathbf{n}} = X_{\mathbf{n}} I\{|X_{\mathbf{n}}| \geq \delta\sqrt{|\mathbf{n}|^{\alpha} \log |\mathbf{n}|}\}.$$

In the following, all objects with primes or multiple primes refer to the respective truncated summands.

Since truncation destroys centering, we obtain, using standard procedures and noting that $EX = 0$,

$$|EX'_{\mathbf{k}}| = |-EX_{\mathbf{k}} I\{|X_{\mathbf{k}}| > b_{\mathbf{k}}\}| \leq E|X|I\{|X_{\mathbf{k}}| > b_{\mathbf{k}}\} \leq \frac{EX^2(\log^+ |X|)^{1-\alpha/2} I\{|X| > b_{\mathbf{k}}\}}{b_{\mathbf{k}}(\log b_{\mathbf{k}})^{1-\alpha/2}}$$

so that

$$\begin{aligned}
|ET'_{\mathbf{n},\mathbf{n}+\mathbf{n}^{\alpha}}| &\leq \sum_{\mathbf{n} \leq \mathbf{k} \leq \mathbf{n}+\mathbf{n}^{\alpha}} \frac{EX^2(\log^+ |X|)^{1-\alpha/2} I\{|X| > b_{\mathbf{k}}\}}{b_{\mathbf{k}}(\log b_{\mathbf{k}})^{1-\alpha/2}} \\
&\leq |\mathbf{n}|^{\alpha} \cdot \frac{EX^2(\log^+ |X|)^{1-\alpha/2} I\{|X| > b_{\mathbf{n}}\}}{b_{\mathbf{n}}(\log b_{\mathbf{n}})^{1-\alpha/2}} \\
&\leq C\sqrt{|\mathbf{n}|^{\alpha}(\log |\mathbf{n}|)^{\alpha}} \cdot EX^2(\log^+ |X|)^{1-\alpha/2} I\{|X| > b_{\mathbf{n}}\} \\
&= o(\sqrt{|\mathbf{n}|^{\alpha} \log |\mathbf{n}|}) \qquad \text{as } \mathbf{n} \to \infty.
\end{aligned} \qquad (4.2)$$

Moreover,

$$\operatorname{Var} X_{\mathbf{n}} \leq EX_{\mathbf{n}}^2 \leq EX^2 = \sigma^2$$

so that

$$\operatorname{Var}(T'_{\mathbf{n},\mathbf{n}+\mathbf{n}^{\alpha}}) \leq |\mathbf{n}|^{\alpha}\sigma^2. \qquad (4.3)$$

An application of the Kolmogorov upper exponential bound (see, e.g., [11], Lemma 8.2.1) with $x = \varepsilon(1-\delta)\sqrt{2\log |\mathbf{n}|}$ and $c_{\mathbf{n}} = 2\delta/x$ (note that $|X'_{\mathbf{k}}| = o(c_n\sqrt{\operatorname{Var}(T_{\mathbf{n}+\mathbf{n}^{\alpha}})})$ for $\mathbf{n} \leq$

$\mathbf{k} \le \mathbf{n} + \mathbf{n}^\alpha$), together with (4.2) and (4.3), now yields

$$P(|T'_{\mathbf{n},\mathbf{n}+\mathbf{n}^\alpha}| > \varepsilon\sqrt{2|\mathbf{n}|^\alpha \log|\mathbf{n}|})$$

$$\le P(|T'_{\mathbf{n},\mathbf{n}+\mathbf{n}^\alpha} - \mathrm{E}T'_{\mathbf{n},\mathbf{n}+\mathbf{n}^\alpha}| > \varepsilon(1-\delta)\sqrt{2|\mathbf{n}|^\alpha \log|\mathbf{n}|})$$

$$\le P\left(|T'_{\mathbf{n},\mathbf{n}+\mathbf{n}^\alpha} - \mathrm{E}T'_{\mathbf{n},\mathbf{n}+\mathbf{n}^\alpha}| > \frac{\varepsilon(1-\delta)}{\sigma}\sqrt{2\,\mathrm{Var}(T'_{\mathbf{n},\mathbf{n}+\mathbf{n}^\alpha})\log|\mathbf{n}|}\right) \qquad (4.4)$$

$$\le \exp\left\{-\frac{2\varepsilon^2(1-\delta)^2}{2\sigma^2}\log|\mathbf{n}|(1-\delta)\right\}$$

$$= |\mathbf{n}|^{-(\varepsilon^2(1-\delta)^3)/\sigma^2}.$$

In order to apply the lower exponential bound (see, e.g., [11], Lemma 8.2.2) we first need a lower bound for the truncated variances:

$$\mathrm{Var}\,X'_{\mathbf{n}} = \mathrm{E}X'^2_{\mathbf{n}} - (\mathrm{E}X'_{\mathbf{n}})^2 = \mathrm{E}X^2 - \mathrm{E}X^2 I\{|X_{\mathbf{n}}| > b_{\mathbf{n}}\} - (\mathrm{E}X'_{\mathbf{n}})^2$$

$$\ge \sigma^2 - 2\mathrm{E}X^2 I\{|X_{\mathbf{n}}| > b_{\mathbf{n}}\} > \sigma^2(1-\delta)$$

for $|\mathbf{n}|$ large, so that

$$\mathrm{Var}(T'_{\mathbf{n},\mathbf{n}+\mathbf{n}^\alpha}) \ge |\mathbf{n}|^\alpha \sigma^2(1-\delta) \qquad \text{for } |\mathbf{n}| \text{ large.} \qquad (4.5)$$

It now follows that for any $\gamma > 0$,

$$P(T'_{\mathbf{n},\mathbf{n}+\mathbf{n}^\alpha} > \varepsilon\sqrt{2|\mathbf{n}|^\alpha \log|\mathbf{n}|})$$

$$\ge P(T'_{\mathbf{n},\mathbf{n}+\mathbf{n}^\alpha} - \mathrm{E}T'_{\mathbf{n},\mathbf{n}+\mathbf{n}^\alpha} > \varepsilon(1+\delta)\sqrt{2|\mathbf{n}|^\alpha \log|\mathbf{n}|})$$

$$\ge P\left(T'_{\mathbf{n},\mathbf{n}+\mathbf{n}^\alpha} - \mathrm{E}T'_{\mathbf{n},\mathbf{n}+\mathbf{n}^\alpha} > \frac{\varepsilon(1+\delta)}{\sigma\sqrt{(1-\delta)}}\sqrt{2\,\mathrm{Var}(T'_{\mathbf{n},\mathbf{n}+\mathbf{n}^\alpha})\log|\mathbf{n}|}\right) \qquad (4.6)$$

$$\ge \exp\left\{-\frac{2\varepsilon^2(1+\delta)^2}{2\sigma^2(1-\delta)}\log|\mathbf{n}|(1+\gamma)\right\}$$

$$= |\mathbf{n}|^{-(\varepsilon^2(1+\delta)^2(1+\gamma))/(\sigma^2(1-\delta))} \qquad \text{for } |\mathbf{n}| \text{ large.}$$

# 5. Proof of Theorem 2.1

We follow the general scheme of [15], although some of the technicalities become more complicated due to the more complicated index set.

## 5.1. Sufficiency – the upper bound

We begin by taking care of the double- and triple-primed contributions, after which we provide a convergent upper Borel–Cantelli sum for the single-primed contribution *over* *a suitably chosen subset* of points in $\mathbf{Z}_+^d$. After this, we apply the first Borel–Cantelli lemma to this subset and then "fill the gaps" in order to include arbitrary windows.

### 5.1.1. $T''_{\mathbf{n},\mathbf{n}+\mathbf{n}^\alpha}$

In this subsection, we establish the fact that

$$\limsup_{\mathbf{n}\to\infty}\frac{|T''_{\mathbf{n},\mathbf{n}+\mathbf{n}^\alpha}|}{\sqrt{|\mathbf{n}|^\alpha\log|\mathbf{n}|}}\le\frac{\delta}{1-\alpha}\qquad\text{a.s.}\tag{5.1}$$

In order for $|T''_{\mathbf{n},\mathbf{n}+\mathbf{n}^\alpha}|$, to surpass the level $\eta\sqrt{|\mathbf{n}|^\alpha\log|\mathbf{n}|}$, it is necessary that at least $N\ge\eta/\delta$ of the $X''$'s are non-zero, which, by stretching the truncation bounds to the extremes, implies that

$$\begin{aligned}
P(|T''_{\mathbf{n},\mathbf{n}+\mathbf{n}^\alpha}|>\eta\sqrt{|\mathbf{n}|^\alpha\log|\mathbf{n}|})&\le\binom{|\mathbf{n}|^\alpha}{N}\left(P(b_{\mathbf{n}}<|X|<\delta\sqrt{(|\mathbf{n}|+|\mathbf{n}|^\alpha)\log(|\mathbf{n}|+|\mathbf{n}|^\alpha)})\right)^N\\
&\le|\mathbf{n}|^{\alpha N}(P(|X|>C|\mathbf{n}|^{\alpha/2}/\log|\mathbf{n}|))^N\\
&\le C|\mathbf{n}|^{\alpha N}\left(\frac{\mathrm{E}|X|^{2/\alpha}(\log^+|X|)^{d-1-1/\alpha}}{(|\mathbf{n}|^{\alpha/2}/\log|\mathbf{n}|)^{2/\alpha}(\log|\mathbf{n}|)^{d-1-1/\alpha}}\right)^N\\
&=C\frac{(\log|\mathbf{n}|)^{N((3/\alpha)+1-d)}}{|\mathbf{n}|^{N(1-\alpha)}}.
\end{aligned}$$

Since the sum of the probabilities converges whenever $N(1-\alpha)>1$, considering that, in addition, $N\delta\ge\eta$, we have shown that

$$\sum_{\mathbf{n}}P(|T''_{\mathbf{n},\mathbf{n}+\mathbf{n}^\alpha}|>\eta\sqrt{|\mathbf{n}|^\alpha\log|\mathbf{n}|})<\infty\qquad\text{for all }\eta>\frac{\delta}{1-\alpha},$$

which establishes (5.1) via the first Borel–Cantelli lemma.

### 5.1.2. $T'''_{\mathbf{n},\mathbf{n}+\mathbf{n}^\alpha}$

Next, we show that

$$\lim_{\mathbf{n}\to\infty}\frac{|T'''_{\mathbf{n},\mathbf{n}+\mathbf{n}^\alpha}|}{\sqrt{|\mathbf{n}|^\alpha\log|\mathbf{n}|}}=0\qquad\text{a.s.}\tag{5.2}$$

This is easier, since, in order for $|T'''_{\mathbf{n},\mathbf{n}+\mathbf{n}^\alpha}|$'s to surpass the level $\eta\sqrt{|\mathbf{n}|^\alpha\log|\mathbf{n}|}$ infinitely often, it is necessary that infinitely many of the $X'''$'s are non-zero. However, via an appeal to the first Borel–Cantelli lemma, the latter event has zero probability since

$$\sum_{\mathbf{n}}P(|X_{\mathbf{n}}|>\eta\sqrt{|\mathbf{n}|^\alpha\log|\mathbf{n}|})=\sum_{\mathbf{n}}P(|X|>\eta\sqrt{|\mathbf{n}|^\alpha\log|\mathbf{n}|})<\infty$$

if and only if (2.1) holds; recall Lemma 3.2.

### 5.1.3. $T'_{\mathbf{n},\mathbf{n}+\mathbf{n}^\alpha}$

As for $T'_{\mathbf{n},\mathbf{n}+\mathbf{n}^\alpha}$, we must resort to subsequences. Set $\lambda_1 = 1$, $\lambda_2 = 2$ and, further,

$$\lambda_i = \left(\frac{i}{\log i}\right)^{1/(1-\alpha)}, \qquad i = 3, 4, \ldots, \quad \text{and} \quad \Lambda = \{\lambda_i, i \geq 1\}.$$

Our attention here is on the subset of points $\mathbf{n} = (n_1, n_2, \ldots, n_d) \in \mathbf{Z}_+^d$ such that $n_k \in \Lambda$, that is, $n_k = \lambda_{i_k}$ for all $k = 1, 2, \ldots, d$, in short $\mathbf{n} \in \boldsymbol{\Lambda}$.

Suppose that $\mathbf{n} \in \boldsymbol{\Lambda}$ and set $i = \prod_{k=1}^d i_k$. This implies, in particular, that $i_k \leq i$ and that $\log i_k \leq \log i$ for all $k$ so that

$$|\mathbf{n}| = \prod_{k=1}^d \lambda_k = \left(\frac{\prod_{k=1}^d i_k}{\prod_{k=1}^d \log i_k}\right)^{1/(1-\alpha)} \geq \frac{i^{1/(1-\alpha)}}{(\log i)^{d/(1-\alpha)}}.$$

With this and (3.3) in mind, the estimate (4.4) over the subset $\boldsymbol{\Lambda}$ now yields

$$\begin{aligned}
&\sum_{\{\mathbf{n} \in \boldsymbol{\Lambda}\}} P(|T'_{\mathbf{n},\mathbf{n}+\mathbf{n}^\alpha}| > \varepsilon\sqrt{2|\mathbf{n}|^\alpha \log |\mathbf{n}|}) \\
&\leq \sum_{\{\mathbf{n} \in \boldsymbol{\Lambda}\}} |\mathbf{n}|^{-(\varepsilon^2(1-\delta)^3)/\sigma^2} \leq \sum_i \sum_{|\prod_{k=1}^d i_k|=i} |\mathbf{n}|^{-(\varepsilon^2(1-\delta)^3)/\sigma^2} \\
&\leq \sum_i d(i)\left(\frac{i^{1/(1-\alpha)}}{(\log i)^{d/(1-\alpha)}}\right)^{-(\varepsilon^2(1-\delta)^3)/(\sigma^2)} \\
&\leq C + \sum_{i \geq i_0} d(i) i^{-(\varepsilon^2((1-\delta)^3 - 2\delta))/(\sigma^2(1-\alpha))} < \infty
\end{aligned} \qquad (5.3)$$

for $\varepsilon > \sigma\sqrt{\frac{1-\alpha}{(1-\delta)^3 - 2\delta}}$ (where $i_0$ was chosen such that $(\log i)^{d(1-\delta)^3} \leq i^\delta$ and $d(i) \leq i^{\delta(\varepsilon^2(1-\delta)^3)/(\sigma^2(1-\alpha))}$ for $i \geq i_0$).

### 5.1.4. Combining the contributions

We first note that an application of the first Borel–Cantelli lemma to (5.3) provides an upper bound for $\limsup T'_{\mathbf{n},\mathbf{n}+\mathbf{n}^\alpha}$ as $\mathbf{n} \to \infty$ through the subset $\boldsymbol{\Lambda}$. More precisely,

$$\limsup_{\substack{\mathbf{n} \to \infty \\ \{\mathbf{n} \in \boldsymbol{\Lambda}\}}} \frac{|T'_{\mathbf{n},\mathbf{n}+\mathbf{n}^\alpha}|}{\sqrt{2|\mathbf{n}|^\alpha \log |\mathbf{n}|}} \leq \sigma\sqrt{\frac{1-\alpha}{(1-\delta)^3 - 2\delta}} \qquad \text{a.s.} \qquad (5.4)$$

Combining this with (5.1) and (5.2) now yields

$$\limsup_{\substack{\mathbf{n} \to \infty \\ \{\mathbf{n} \in \boldsymbol{\Lambda}\}}} \frac{|T_{\mathbf{n},\mathbf{n}+\mathbf{n}^\alpha}|}{\sqrt{2|\mathbf{n}|^\alpha \log |\mathbf{n}|}} \leq \sigma\sqrt{\frac{1-\alpha}{(1-\delta)^3 - 2\delta}} + \frac{\delta}{1-\alpha} \qquad \text{a.s.,}$$

which, due to the arbitrary nature of $\delta$, tells us that

$$\limsup_{\substack{\mathbf{n}\to\infty \\ \{\mathbf{n}\in\boldsymbol{\Lambda}\}}} \frac{|T_{\mathbf{n},\mathbf{n}+\mathbf{n}^\alpha}|}{\sqrt{2|\mathbf{n}|^\alpha \log|\mathbf{n}|}} \le \sigma\sqrt{1-\alpha} \qquad \text{a.s.} \tag{5.5}$$

### 5.1.5. *Filling the gaps*

Thus far, we have shown that the lim sup for a special *subsequence* is as desired. We must now show that the lim sup remains the same for an arbitrary sequence. For the usual LIL, this is typically done by studying the gaps between subsequence points with the aid of the Lévy inequalities. Here, however, we must proceed differently.

Since the procedure is the same in all dimensions, we restrict ourselves to carrying out the details for the case $d=2$. Since, in the following, we shall use the letters $m$ and $n$ for $x$- and $y$-coordinates, respectively, we set

$$m_j = n_j = \lambda_j, \qquad j \ge 1,$$

that is, the points $(m_j, n_k)$ for $j, k = 1, 2, \dots$ are the Southwest corners of the windows we have considered thus far.

The first step is to show that the selected windows overlap, that is, that they cover all of $\mathbf{Z}_+^d$. For this purpose, it suffices to consider squares. We thus wish to show that

$$m_i + m_i^\alpha > m_{i+1} \qquad \text{for all } i,$$

that is, that

$$\left(\frac{i}{\log i}\right)^{1/(1-\alpha)} + \left(\frac{i}{\log i}\right)^{\alpha/(1-\alpha)} > \left(\frac{i+1}{\log(i+1)}\right)^{1/(1-\alpha)} \qquad \text{for all } i \ge 3, \tag{5.6}$$

in other words, that the Northeast end-point of one square overlaps the Southwest end-point of the following Northeast-square. This, however, follows from the fact that

$$\frac{(i/(\log i))^{1/(1-\alpha)} + (i/(\log i))^{\alpha/(1-\alpha)}}{((i+1)/(\log(i+1)))^{1/(1-\alpha)}}$$

$$= \left(\frac{i\log(i+1)}{(i+1)\log i}\right)^{1/(1-\alpha)} \cdot \left(1 + \frac{\log i}{i}\right)$$

$$= \left(\frac{\log(i+1)}{\log i}\right)^{1/(1-\alpha)} \cdot \left(1 - \frac{1}{i+1}\right)^{1/(1-\alpha)} \cdot \left(1 + \frac{\log i}{i}\right)^{1/(1-\alpha)}$$

$$\ge \left(1 - \frac{1}{i+1}\right)^{1/(1-\alpha)} \cdot \left(1 + \frac{\log i}{i}\right)^{1/(1-\alpha)}$$

$$\ge \left(1 + \frac{1}{i}\left(\log i - \frac{i}{i+1} - \frac{\log i}{i(i+1)}\right)\right)^{1/(1-\alpha)} > 1 \qquad \text{for } i \ge 3.$$

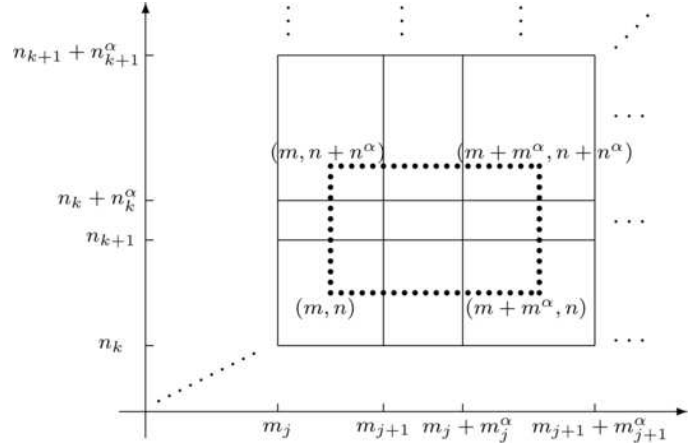

**Figure 1.** A dotted arbitrary window.

Next, we select an arbitrary window,

$$T_{((m,n),(m+m^\alpha, n+n^\alpha))}.$$

Since, trivially, we have $n_k + n_k^\alpha \leq n_{k+1} + n_{k+1}^\alpha$ (as well as $m_j + m_j^\alpha \leq m_{j+1} + m_{j+1}^\alpha$), it follows that an arbitrary window is always contained in the union of (at most) four selected windows as depicted in Figure 1.

The program for this subsubsection is (essentially) to show that the discrepancy between an arbitrary window and the original ones is asymptotically negligible.

From (5.4), we recall that

$$\limsup_{\substack{\mathbf{n}\to\infty \\ \{\mathbf{n}:|\mathbf{n}|=\lambda_i\}}} \frac{|T'_{\mathbf{n},\mathbf{n}+\mathbf{n}^\alpha}|}{\sqrt{2|\mathbf{n}|^\alpha \log|\mathbf{n}|}} \leq \sigma \sqrt{\frac{1-\alpha}{(1-\delta)^3 - 2\delta}} \qquad \text{a.s.}$$

We wish to show that the same relation holds for the full subsequence, that is, that

$$\limsup_{\mathbf{n}\to\infty} \frac{|T'_{\mathbf{n},\mathbf{n}+\mathbf{n}^\alpha}|}{\sqrt{2|\mathbf{n}|^\alpha \log|\mathbf{n}|}} \leq \sigma \sqrt{\frac{1-\alpha}{(1-\delta)^3 - 2\delta}} \qquad \text{a.s.,}$$

which, remembering that we are restricting ourselves to the case $d = 2$, transforms into

$$\limsup_{j,k\to\infty} \frac{\sum_{i_1=m_j}^{m_j+m_j^\alpha} \sum_{i_2=n_k}^{n_k+n_k^\alpha} X'_{i_1,i_2}}{\sqrt{2m_j^\alpha n_k^\alpha \log(m_j n_k)}} \leq \sigma \sqrt{\frac{1-\alpha}{(1-\delta)^3 - 2\delta}} \qquad \text{a.s.}$$

by showing that

$$\limsup_{j,k\to\infty} \max_{\substack{m_j<m\le m_{j+1}\\ n_k<n\le n_{j+1}}} \left| \frac{\sum_{i_1=m}^{m+m^\alpha}\sum_{i_2=n}^{n+n^\alpha} X'_{i_1,i_2}}{\sqrt{2m^\alpha n^\alpha \log(mn)}} - \frac{\sum_{i_1=m_j}^{m_j+m_j^\alpha}\sum_{i_2=n_k}^{n_k+n_k^\alpha} X'_{i_1,i_2}}{\sqrt{2m_j^\alpha n_k^\alpha \log(m_j n_k)}} \right| = 0 \qquad \text{a.s.}$$

However, since $n_{k+1}/n_k \to 1$ as $k\to\infty$ (and $m_{j+1}/m_j \to 1$ as $j\to\infty$) it suffices to show, say, that

$$\limsup_{j,k\to\infty} \max_{\substack{m_j<m\le m_{j+1}\\ n_k<n\le n_{k+1}}} \frac{|\sum_{i_1=m}^{m+m^\alpha}\sum_{i_2=n}^{n+n^\alpha} X'_{i_1,i_2} - \sum_{i_1=m_j}^{m_j+m_j^\alpha}\sum_{i_2=n_k}^{n_k+n_k^\alpha} X'_{i_1,i_2}|}{\sqrt{m_j^\alpha n_k^\alpha \log(m_j n_k)}} = 0 \qquad \text{a.s. (5.7)}$$

By combining this with our previous results concerning $T'''_{\mathbf{n},\mathbf{n}+\mathbf{n}^\alpha}$ and $T''''_{\mathbf{n},\mathbf{n}+\mathbf{n}^\alpha}$, as in the previous subsubsection, we are then in the position to conclude that relation (5.5) holds for the full sequence, that is,

$$\limsup_{\mathbf{n}\to\infty} \frac{|T_{\mathbf{n},\mathbf{n}+\mathbf{n}^\alpha}|}{\sqrt{2|\mathbf{n}|^\alpha \log|\mathbf{n}|}} \le \sigma\sqrt{1-\alpha} \qquad \text{a.s.},$$

as desired (since now $d=2$).

In order to pursue our task, we separate the index set into 3 pieces (since $d=2$) depending on whether the arbitrary window is located in "the center" or "near" one of the coordinate axes (for a similar discussion, cf. [8], Section 4).

*The center; $j,k \ge M$, $M$ large*

**Proof of 5.7.** Let $D_{m,n}$ denote the random variable in the numerator of (5.7), that is, set

$$\begin{aligned}
D_{m,n} &= \sum_{i_1=m+1}^{m+m^\alpha}\sum_{i_2=n+1}^{n+n^\alpha} X'_{i_1,i_2} - \sum_{i_1=m_j+1}^{m_j+m_j^\alpha}\sum_{i_2=n_k+1}^{n_k+n_k^\alpha} X'_{i_1,i_2} \\
&= \sum_{i_1=m_j+m_j^\alpha+1}^{m+m^\alpha}\sum_{i_2=n+1}^{n+n^\alpha} X'_{i_1,i_2} + \sum_{i_1=m+1}^{m_j+m_j^\alpha}\sum_{i_2=n_k+n_k^\alpha+1}^{n+n^\alpha} X'_{i_1,i_2} \qquad (5.8) \\
&\quad - \sum_{i_1=m_j+1}^{m_j+m_j^\alpha}\sum_{i_2=n_k+1}^{n} X'_{i_1,i_2} - \sum_{i_1=m_j+1}^{m}\sum_{i_2=n+1}^{n_k+n_k^\alpha} X'_{i_1,i_2}.
\end{aligned}$$

The corresponding random variables are located in the shaded area in Figure 2.

As for the truncated means and variances, we recall that

$$\mathrm{E}X'_{\mathbf{k}} = o(\sqrt{\log(|\mathbf{k}|)/|\mathbf{k}|^\alpha}) \qquad \text{as } \mathbf{k}\to\infty, \quad \text{and that} \quad \mathrm{Var}\, X'_{\mathbf{k}} \le \mathrm{E}(X'_{\mathbf{k}})^2 \le \sigma^2,$$

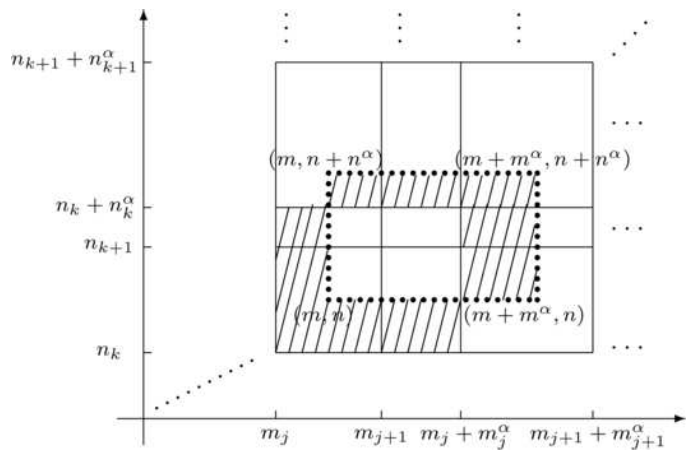

**Figure 2.** The shaded difference.

which, in the present case, means that

$$\mathrm{E}X'_{i_1,i_2} = o\big(\sqrt{\log(m_j n_k)/(m_j^\alpha n_k^\alpha)}\big) \qquad \text{as } j,k \to \infty, \quad \text{and that} \quad \mathrm{Var}\, X'_{i_1,i_2} \leq \sigma^2, \quad (5.9)$$

uniformly in $i_1 \geq m_j$ and $i_2 \geq n_k$.

In order to compute the mean and variance of $D_{m,n}$, we need an estimate of the number of summands involved in (5.8).

As a preliminary, we note that, via the mean value theorem,

$$m_{j+1} - m_j = \left(\frac{j+1}{\log(j+1)}\right)^{1/(1-\alpha)} - \left(\frac{j}{\log j}\right)^{1/(1-\alpha)}$$

$$\sim \frac{1}{1-\alpha} \frac{j^{\alpha/(1-\alpha)}}{(\log j)^{1/(1-\alpha)}} \qquad \text{as } j \to \infty,$$

$$m_{j+1}^\alpha - m_j^\alpha = \left(\frac{j+1}{\log(j+1)}\right)^{\alpha/(1-\alpha)} - \left(\frac{j}{\log j}\right)^{\alpha/(1-\alpha)}$$

$$\sim \frac{\alpha}{1-\alpha} \frac{j^{(2\alpha-1)/(1-\alpha)}}{(\log j)^{\alpha/(1-\alpha)}} \qquad \text{as } j \to \infty.$$

Since $2\alpha - 1 < \alpha$, we also note that

$$m_{j+1}^\alpha - m_j^\alpha = o(m_{j+1} - m_j) \qquad \text{as } j \to \infty.$$

Analogous relations obviously also hold for $n_{k+1} - n_k$ and $n_{k+1}^\alpha - n_k^\alpha$.

Using all of this, we find that the number of summands in (5.8) equals

$$
\begin{aligned}
\mathrm{Card}(D_{m,n}) &= (m + m^\alpha - m_j - m_j^\alpha)n^\alpha + (m_j + m_j^\alpha - m)(n + n^\alpha - n_k - n_k^\alpha) \\
&\quad + m_j^\alpha(n - n_k) + (m - m_j)(n_k + n_k^\alpha - n) \\
&\leq (m_{j+1} + m_{j+1}^\alpha - m_j - m_j^\alpha)n_{k+1} + m_j^\alpha(n_{k+1} + n_{k+1}^\alpha - n_k - n_k^\alpha) \\
&\quad + m_j^\alpha(n_{k+1} - n_k) + (m_{j+1} - m_j)n_k^\alpha \\
&\sim C \frac{j^{\alpha/(1-\alpha)}}{(\log j)^{1/(1-\alpha)}} \left(\frac{k}{\log k}\right)^{\alpha/(1-\alpha)} + \frac{k^{\alpha(1-\alpha)}}{(\log k)^{1/(1-\alpha)}} \left(\frac{j}{\log j}\right)^{\alpha/(1-\alpha)} \\
&= C m_j^\alpha n_k^\alpha \left(\frac{1}{\log j} + \frac{1}{\log k}\right) \qquad \text{as } j, k \to \infty.
\end{aligned}
$$

Combining this with the estimate for the truncated expectations and variances in (5.9), we conclude that, for $m_j \leq m \leq m_{j+1}$, $n_k \leq n \leq n_{k+1}$ and $j, k \to \infty$,

$$
\mathrm{E}(D_{m,n}) = o\left(\sqrt{m_j^\alpha n_k^\alpha \log(m_j n_k) \left(\frac{1}{\log j} + \frac{1}{\log k}\right)}\right)
$$

and

$$
\mathrm{Var}(D_{m,n}) \leq C m_j^\alpha n_k^\alpha \left(\frac{1}{\log j} + \frac{1}{\log k}\right)\sigma^2.
$$

Now, let $\eta > 0$ be arbitrarily small. Using the upper exponential inequalities, we obtain, for $j, k$ large,

$$
\begin{aligned}
&P\left(\frac{|D_{m,n}|}{\sqrt{m_j^\alpha n_k^\alpha \log(m_j n_k)}} > 2\eta\right) \\
&= P\left(|D_{m,n} - \mathrm{E}D_{m,n}| > \eta \sqrt{\frac{m_j^\alpha n_k^\alpha \log(m_j n_k)}{\mathrm{Var}\, D_{m,n}}} \cdot \sqrt{\mathrm{Var}\, D_{m,n}}\right) \\
&\leq P\left(|D_{m,n} - \mathrm{E}D_{m,n}| > C\eta \sqrt{\frac{\log j \log k \log(m_j n_k)}{\log j + \log k}} \cdot \sqrt{\mathrm{Var}\, D_{m,n}}\right) \\
&\sim P\left(|D_{m,n} - \mathrm{E}D_{m,n}| > C\eta \sqrt{\log j \log k} \sqrt{\mathrm{Var}\, D_{m,n}}\right) \\
&\leq \exp\left\{-\frac{1}{2} C\eta^2 \log j \log k\right\},
\end{aligned}
$$

independent of $(n, m) \in [m_j, m_{j+1}] \times [n_k, n_{k+1}]$. Next, let $M$ large be given and let $j, k > e^{2M}$. Then,

$$
\log j \log k \geq M(\log j + \log k),
$$

which implies that

$$\exp\{-\tfrac{1}{2}C\eta^2 \log j \log k\} \le \exp\{-\tfrac{1}{2}C\eta^2 M(\log j + \log k)\} = (jk)^{-C\eta^2 M}$$

and, hence, that

$$P\left(\max_{\substack{m_j < m \le m_{j+1} \\ n_k < n \le n_{k+1}}} \frac{|D_{m,n}|}{\sqrt{m_j^\alpha n_k^\alpha \log(m_j n_k)}} > 2\eta\right)$$

$$\le (m_{j+1} - m_j)(n_{k+1} - n_k) \max_{\substack{m_j < m \le m_{j+1} \\ n_k < n \le n_{j+1}}} P\left(\frac{|D_{m,n}|}{\sqrt{m_j^\alpha n_k^\alpha \log(m_j n_k)}} > 2\eta\right)$$

$$\le C \frac{j^{\alpha/(1-\alpha)}}{(\log j)^{1/(1-\alpha)}} \frac{k^{\alpha/(1-\alpha)}}{(\log k)^{1/(1-\alpha)}} (kj)^{-C\eta^2 M}$$

$$\le C(jk)^{-C\eta^2 M},$$

so that, by choosing $M$ sufficiently large that $C\eta^2 M > 1$, we may finally conclude that

$$\sum_{j,k} P\left(\max_{\substack{m_j < m \le m_{j+1} \\ n_k < n \le n_{k+1}}} \frac{|D_{m,n}|}{\sqrt{m_j^\alpha n_k^\alpha \log(m_j n_k)}} > \eta\right) < \infty \qquad \text{for any } \eta > 0,$$

which, in turn, verifies (5.7) – for an arbitrary "central" window. □

*The case $k \le M$, $M$ large*

For the "boundary" cases, we consider a denser subsequence, namely

$$A_{j,k} = \left\{(j,k) \in \mathbf{Z}_+^2 : m_j = \left(\frac{j}{\log j}\right)^{1/(1-\alpha)}, \ j \ge 3, \text{ and } k = 1, 2, \dots, M\right\}.$$

A consequence of this is that additional windows are involved so that we must first convince ourselves that the upper bound of the lim sup of the thus chosen subset remains the same. Since, as we have seen, the double- and triple-primed contributions do not contribute, it follows that it suffices to investigate $T'_{\mathbf{n}, \mathbf{n}+\mathbf{n}^\alpha}$. In fact, borrowing from (4.4), we have, for $\mathbf{n} = (m_j, k)$ large,

$$P(|T'_{\mathbf{n},\mathbf{n}+\mathbf{n}^\alpha}| > \varepsilon\sqrt{2\mathbf{n}^\alpha \log|\mathbf{n}|}) \le |\mathbf{n}|^{-(\varepsilon^2(1-\delta)^3)/\sigma^2} = (m_j k)^{-(\varepsilon^2(1-\delta)^3)/\sigma^2}$$

so that

$$\sum_{\{\mathbf{n}=(j,k)\in A_{j,k}\}} P(|T'_{\mathbf{n},\mathbf{n}+\mathbf{n}^\alpha}| > \epsilon\sqrt{2\mathbf{n}^\alpha \log|\mathbf{n}|})$$

$$\le \sum_{j=3}^{\infty} \sum_{k=1}^{M} \left(\left(\frac{j}{\log j}\right)^{1/(1-\alpha)} k\right)^{-(\varepsilon^2(1-\delta)^3)/\sigma^2}$$

$$\leq \tilde{M} \sum_{j=3}^{\infty} \left( \frac{j}{\log j} \right)^{-(\varepsilon^2 (1-\delta)^3)/(\sigma^2(1-\alpha))} < \infty \qquad \text{for } \varepsilon > \sqrt{\frac{\sigma^2(1-\alpha)}{(1-\delta)^3}},$$

which is the same as

$$\limsup_{\substack{j \to \infty \\ \{(j,k) \in A_{j,k}\}}} \frac{\sum_{i_1=n_j}^{n_j+n_j^\alpha} \sum_{i_2=k}^{k+k^\alpha} X'_{i_1,i_2}}{\sqrt{2n_j^\alpha k^\alpha \log(n_j k)}} \leq \sigma \sqrt{\frac{1-\alpha}{(1-\delta)^3}} \qquad \text{a.s.} \tag{5.10}$$

To complete the proof for this case, it remains to establish the analog of (5.7), that is, that

$$\limsup_{j \to \infty} \max_{\substack{m_j < m \leq m_{j+1} \\ n \leq M}} \frac{|\sum_{i_1=m}^{m+m^\alpha} \sum_{i_2=n}^{n+n^\alpha} X'_{i_1,i_2} - \sum_{i_1=m_j}^{m_j+m_j^\alpha} \sum_{i_2=n}^{n+n^\alpha} X'_{i_1,i_2}|}{\sqrt{m_j^\alpha n^\alpha \log(m_j n)}} = 0 \qquad \text{a.s. (5.11)}$$

**Proof of (5.11).** With $D_{m,n}$ as the random variable in the numerator of (5.11), we have

$$D_{m,n} = \sum_{i_1=m+1}^{m+m^\alpha} \sum_{i_2=n+1}^{n+n^\alpha} X'_{i_1,i_2} - \sum_{i_1=m_j+1}^{m_j+m_j^\alpha} \sum_{i_2=n}^{n+n^\alpha} X'_{i_1,i_2} = \sum_{i_1=m_j}^{m} \sum_{i_2=n+1}^{n+n^\alpha} X'_{i_1,i_2} + \sum_{i_1=m_j+m_j^\alpha}^{m+m^\alpha} \sum_{i_2=n+1}^{n+n^\alpha} X'_{i_1,i_2},$$

which corresponds to the shaded area in Figure 3.

Continuing as before, we find that

$$\text{Card}(D_{m,n}) = (m + m^\alpha - m_j - m_j^\alpha)n^\alpha + (m - m_j)n^\alpha$$

$$\leq (m_{j+1} + m_{j+1}^\alpha - m_j - m_j^\alpha)n^\alpha + (m_{j+1} - m_j)n^\alpha$$

$$\sim C \frac{j^{\alpha/(1-\alpha)}}{(\log j)^{1/(1-\alpha)}} n^\alpha = C \frac{m_j^\alpha n^\alpha}{\log j}.$$

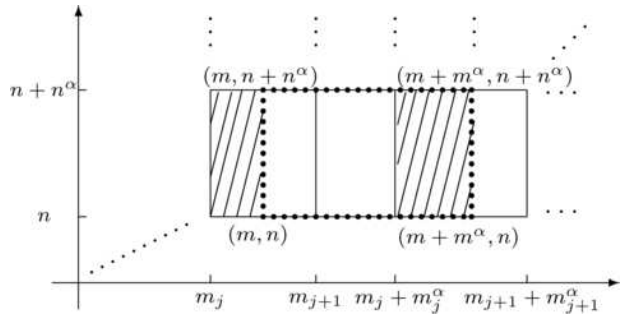

**Figure 3.** The shaded difference.

Combining this with the estimate for the truncated expectations (which, in this case, corresponds to $\mathrm{E}X_{i_1,n} = o(\sqrt{\log(m_j n)/m_j^\alpha n^\alpha})$) and variances in (5.9), we conclude that, for $m_j \le m \le m_{j+1}$, $n \le M$ and $j, n \to \infty$,

$$\mathrm{E}(D_{m,n}) = o\left(\sqrt{\frac{m_j^\alpha n^\alpha \log(m_j n)}{\log j}}\right) \quad \text{and} \quad \mathrm{Var}(D_{m,n}) \le C\frac{m_j^\alpha n^\alpha}{\log j}\sigma^2.$$

The exponential inequalities therefore yield

$$P\left(\frac{|D_{m,n}|}{\sqrt{m_j^\alpha n^\alpha \log(m_j n)}} > 2\eta\right) = P\left(|D_{m,n} - \mathrm{E}D_{m,n}| > \eta\sqrt{\frac{m_j^\alpha n^\alpha \log(m_j n)}{\mathrm{Var}\,D_{m,n}}} \cdot \sqrt{\mathrm{Var}\,D_{m,n}}\right)$$

$$\le P(|D_{m,n} - \mathrm{E}D_{m,n}| > C\eta\sqrt{\log j \log(m_j n)} \cdot \sqrt{\mathrm{Var}\,D_{m,n}})$$

$$\sim P(|D_{m,n} - \mathrm{E}D_{m,n}| > C\eta \log j \sqrt{\log n}\sqrt{\mathrm{Var}\,D_{m,n}})$$

$$\le \exp\left\{-\frac{1}{2}C\eta^2(\log j)^2 \log n\right\} \le \exp\left\{-\frac{1}{2}C\eta^2(\log j)^2\right\}$$

so that

$$P\left(\max_{\substack{m_j < m \le m_{j+1} \\ n \le M}} \frac{|D_{m,n}|}{\sqrt{m_j^\alpha n^\alpha \log(m_j n)}} > 2\eta\right)$$

$$\le (m_{j+1} - m_j)M \max_{\substack{m_j < m \le m_{j+1} \\ n \le M}} P\left(\frac{|D_{m,n}|}{\sqrt{m_j^\alpha n^\alpha \log(m_j n)}} > 2\eta\right)$$

$$\le C\frac{j^{\alpha/(1-\alpha)}}{(\log j)^{1/(1-\alpha)}}\exp\left\{-\frac{1}{2}C\eta^2(\log j)^2\right\}$$

$$\le \exp\left\{-\frac{1}{2}C\eta^2(\log j)^2 + \frac{1}{1-\alpha}\log j\right\}$$

$$\le \exp\left\{-\frac{1}{2}C\eta^2(\log j)^2\right\},$$

from which we finally conclude that

$$\sum_j P\left(\max_{\substack{m_j < m \le m_{j+1} \\ n \le M}} \frac{|D_{m,n}|}{\sqrt{m_j^\alpha n^\alpha \log(m_j n)}} > \eta\right) < \infty \qquad \text{for any } \eta > 0,$$

which, in turn, verifies (5.11).

*The case $j \le M$, $M$ large*

This part clearly follows in the same way as the previous one by interchanging the roles of $j$ and $k$.                                                                                                        □

This, finally, concludes the proof of the upper bound

$$\limsup_{\mathbf{n}\to\infty} \frac{|T_{\mathbf{n},\mathbf{n}+\mathbf{n}^\alpha}|}{\sqrt{|\mathbf{n}|^\alpha \log |\mathbf{n}|}} \leq \sigma\sqrt{1-\alpha} \qquad \text{a.s.} \tag{5.12}$$

## 5.2. Sufficiency – the lower bound

We first derive a divergent Borel–Cantelli sum for the single-primed contributions restricted to the subsets of windows based on the subsequence $\lambda_1 = 1$, $\lambda_2 = 2$ and $\lambda_i = i^{1/(1-\alpha)}$, $i \geq 3$. More precisely, in dimension 2, the Southwest coordinates are $(i^{1/(1-\alpha)}, k^{1/(1-\alpha)})$, the horizontal widths are $i^{\alpha/(1-\alpha)}$ and the vertical widths are $k^{\alpha/(1-\alpha)}$, $i,k \geq 1$. (See Figure 4.) In order to apply the second Borel–Cantelli lemma to this subset, we then show that this provides an *independent subset* of windows. Combining this with (5.1) and (5.2) provides a lower bound for the subset and, further, for the whole set.

### 5.2.1. $T'_{\mathbf{n},\mathbf{n}+\mathbf{n}^\alpha}$

Let $A$ denote the set of Southwest coordinates involved, that is, set

$$A = \{\mathbf{n} \in \mathbf{Z}_+^2 : \mathbf{n} = (i_1^{1/(1-\alpha)}, i_2^{1/(1-\alpha)}),\ i_1, i_2 \geq 1\}.$$

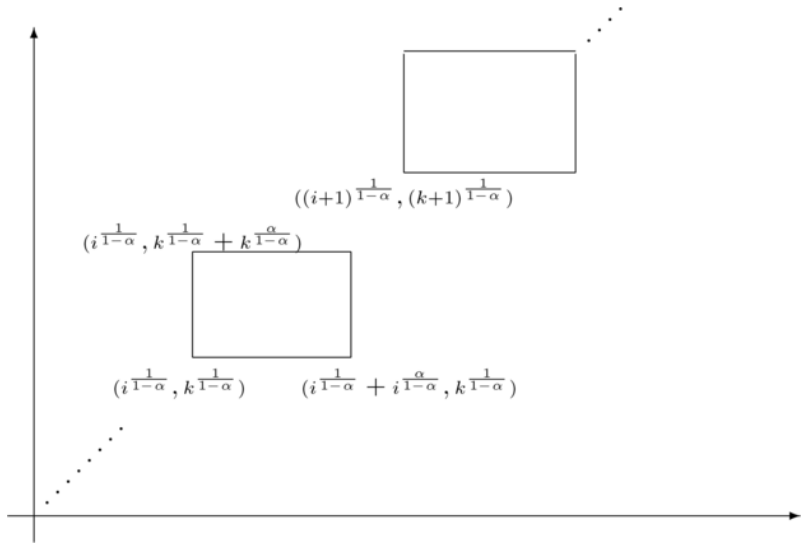

**Figure 4.** Independence of windows.

An application of (4.6) now tells us that

$$
\begin{aligned}
\sum_{\{\mathbf{n} \in A\}} P(T'_{\mathbf{n}, \mathbf{n}+\mathbf{n}^\alpha} > \varepsilon \sqrt{|\mathbf{n}|^\alpha \log |\mathbf{n}|}) &\geq \sum_{\{\mathbf{n} \in A\}} |\mathbf{n}|^{-(\varepsilon^2(1+\delta)^2(1+\gamma))/(\sigma^2(1-\delta))} \\
&= \sum_i \sum_{\{i_1, i_2 : i_1 \cdot i_2 = i\}} |\mathbf{n}|^{-(\varepsilon^2(1+\delta)^2(1+\gamma))/(\sigma^2(1-\delta))} \\
&= \sum_i d(i) i^{-(\varepsilon^2(1+\delta)^2(1+\gamma))/(\sigma^2(1-\alpha)(1-\delta))} = \infty
\end{aligned}
\tag{5.13}
$$

for all $\varepsilon < \sigma \sqrt{\frac{(1-\alpha)(1-\delta)}{(1+\delta)^2(1+\gamma)}}$.

### 5.2.2. Independence

In order to prove that the selected windows are disjoint, it suffices to check that, throughout, each coordinate with index $i+1$ is larger than the corresponding coordinate of $\mathbf{n} + \mathbf{n}^\alpha$ with index $i$. This means that we must show that

$$
i^{1/(1-\alpha)} + i^{\alpha/(1-\alpha)} < (i+1)^{1/(1-\alpha)} \qquad \text{for all } i.
\tag{5.14}
$$

However, this follows from the fact that

$$
\frac{i^{1/(1-\alpha)} + i^{\alpha/(1-\alpha)}}{(i+1)^{1/(1-\alpha)}} = \left(\frac{i}{i+1}\right)^{1/(1-\alpha)} \cdot \left(1 + \frac{1}{i}\right) = \left(\frac{i}{i+1}\right)^{\alpha/(1-\alpha)} < 1.
$$

### 5.2.3. Combining the contributions

We have just shown that the selected subset of windows is disjoint and, hence, that the events considered in (5.13) are independent. An application of the second Borel–Cantelli lemma therefore tells us that

$$
\limsup_{\substack{\mathbf{n} \to \infty \\ \{\mathbf{n}: \mathbf{n} \in A\}}} \frac{T'_{\mathbf{n}, \mathbf{n}+\mathbf{n}^\alpha}}{\sqrt{|\mathbf{n}|^\alpha \log |\mathbf{n}|}} \geq \sigma \sqrt{\frac{(1-\alpha)(1-\delta)}{(1+\delta)^2(1+\gamma)}}.
$$

Combining this with (5.1) and (5.2), it follows that

$$
\limsup_{\substack{\mathbf{n} \to \infty \\ \{\mathbf{n}: \mathbf{n} \in A\}}} \frac{T_{\mathbf{n}, \mathbf{n}+\mathbf{n}^\alpha}}{\sqrt{|\mathbf{n}|^\alpha \log |\mathbf{n}|}} \geq \sigma \sqrt{\frac{(1-\alpha)(1-\delta)}{(1+\delta)^2(1+\gamma)}} - \frac{\delta}{1-\alpha} \qquad \text{a.s.}
$$

Therefore, due to the arbitrary nature of $\delta$ and $\gamma$, it follows that

$$
\limsup_{\substack{\mathbf{n} \to \infty \\ \{\mathbf{n}: \mathbf{n} \in A\}}} \frac{T_{\mathbf{n}, \mathbf{n}+\mathbf{n}^\alpha}}{\sqrt{|\mathbf{n}|^\alpha \log |\mathbf{n}|}} \geq \sigma \sqrt{1-\alpha} \qquad \text{a.s.}
\tag{5.15}
$$

and, since the overall lim sup is obviously at least as large as the selected one, that

$$\limsup_{\mathbf{n}\to\infty} \frac{T_{\mathbf{n},\mathbf{n}+\mathbf{n}^\alpha}}{\sqrt{|\mathbf{n}|^\alpha \log|\mathbf{n}|}} \geq \sigma\sqrt{1-\alpha} \qquad \text{a.s.} \tag{5.16}$$

*5.2.4. Final step*

The proof of the sufficiency is completed by combining (5.12) and (5.16).

## 5.3. Necessity

If (2.3) holds, then, by the zero–one law, the probability that the lim sup is finite is 0 or 1. Hence, being positive, it equals 1. Consequently (cf. [15], page 438 or [17, 18]),

$$\limsup_{n\to\infty} \frac{|X_\mathbf{n}|}{\sqrt{|\mathbf{n}|^\alpha \log|\mathbf{n}|}} < \infty \qquad \text{a.s.,}$$

from which it follows, via the second Borel–Cantelli lemma and the i.i.d. assumption, that

$$\infty > \sum_\mathbf{n} P(|X_\mathbf{n}| > \sqrt{|\mathbf{n}|^\alpha \log|\mathbf{n}|})$$
$$= \sum_\mathbf{n} P(|X| > \sqrt{|\mathbf{n}|^\alpha \log|\mathbf{n}|}).$$

This verifies (2.1) in view of Lemma 3.2 (with $\beta = 1$).

An application of the sufficiency part finally tells us that (2.2) holds with $\sigma^2 = \text{Var}\, X$.

# 6. An LSL for subsequences

As we have seen, it follows from (5.1) and (5.2) that the values of the extreme limit points are determined by the behaviour of $T'_{\mathbf{n},\mathbf{n}+\mathbf{n}^\alpha}$. In this section, we shall exploit this fact further and prove an LSL for subsequences (paralleling [10], where this was done for the classical LIL; cf. also [11], Section 8.5).

To this end, we replace the set $\Lambda$ in the computation of the upper bound by the set

$$\Lambda^* = \{\lambda_i^* = \lambda_i^\beta = i^{\beta/(1-\alpha)}, i \geq 1\} \tag{6.1}$$

and, consequently, $\mathbf{\Lambda}$ by $\mathbf{\Lambda}^*$.

Now, if all coordinates of $\mathbf{n}$ belong to $\Lambda^*$ (in short, $\mathbf{n} \in \mathbf{\Lambda}^*$), then so does $|\mathbf{n}|$, that is, the "size" of the points in $\mathbf{Z}_+^d$. This means that $\sum_{\{\mathbf{n}\in\mathbf{\Lambda}^*\}} = \sum_{\{\mathbf{n}:|\mathbf{n}|\in\mathbf{\Lambda}^*\}}$, that is, summation runs over a sequence of hyperbolas approaching infinity.

With this in mind, an inspection of the computation of the upper bound of $T'_{\mathbf{n},\mathbf{n}+\mathbf{n}^\alpha}$, together with an application of Lemma 3.3 with $\kappa = \beta/(1-\alpha)$, $\theta = \varepsilon^2(1-\delta)^3/\sigma^2$ and $\eta = 0$, reveals that formula (5.3) becomes

$$\sum_{\{\mathbf{n}\in\boldsymbol{\Lambda}^*\}} P(|T'_{\mathbf{n},\mathbf{n}+\mathbf{n}^\alpha}| > \varepsilon\sqrt{|\mathbf{n}|^\alpha\log|\mathbf{n}|})$$

$$= \sum_i \sum_{\{\mathbf{n}:|\mathbf{n}|=i^{\beta/(1-\alpha)}\}} P(|T'_{\mathbf{n},\mathbf{n}+\mathbf{n}^\alpha}| > \varepsilon\sqrt{|\mathbf{n}|^\alpha\log|\mathbf{n}|})$$

$$\leq \sum_i \sum_{\{\mathbf{n}:|\mathbf{n}|=i^{\beta/(1-\alpha)}\}} |\mathbf{n}|^{-(\varepsilon^2(1-\delta)^3)/\sigma^2} < \infty \qquad (6.2)$$

for all $\varepsilon > \sigma\sqrt{\frac{1-\alpha}{\beta(1-\delta)^3}}$. Moreover, with $i$ raised to the various powers in the computations of the lower bound replaced by $i$ raised to $\beta$ times the same powers and the index set $A$ replaced by

$$A^* = \{\mathbf{n}\in\mathbf{Z}_+^2 : \mathbf{n} = (i^{\beta/(1-\alpha)}, k^{\beta/(1-\alpha)}), \ i\geq 1\}, \qquad (6.3)$$

(5.13) transforms into

$$\sum_{\{\mathbf{n}\in A^*\}} P(T'_{\mathbf{n},\mathbf{n}+\mathbf{n}^\alpha} > \varepsilon\sqrt{|\mathbf{n}|^\alpha\log|\mathbf{n}|}) \geq \sum_{\{\mathbf{n}\in A^*\}} |\mathbf{n}|^{-(\varepsilon^2\beta(1+\delta)^2(1+\gamma))/(\sigma^2(1-\delta))}$$

$$= \sum_i d(i) i^{-(\varepsilon^2\beta(1+\delta)^2(1+\gamma))/(\sigma^2(1-\alpha)(1-\delta))} = \infty \qquad (6.4)$$

for all $\varepsilon < \sigma\sqrt{\frac{(1-\alpha)(1-\delta)}{\beta(1+\delta)^2(1+\gamma)}}$.

By combining these estimates with (5.1) and (5.2) as in the proof of Theorem 2.1, we obtain the following *LSL for subsequences*.

**Theorem 6.1.** *Suppose that $\{X_{\mathbf{k}}, \mathbf{k}\in\mathbf{Z}_+^d\}$ are i.i.d. random variables with mean 0 and finite variance $\sigma^2$, set $S_{\mathbf{n}} = \sum_{\mathbf{k}\leq\mathbf{n}} X_{\mathbf{k}}$, $\mathbf{n}\in\mathbf{Z}_+^d$, and let $\boldsymbol{\Lambda}^*$ be as defined in (6.1). If*

$$\mathrm{E} X^{2/\alpha}(\log^+|X|)^{d-1-1/\alpha} < \infty, \qquad (6.5)$$

*where $0 < \alpha < 1$, then, for $\beta > 1$,*

$$\limsup_{\substack{\mathbf{n}\to\infty \\ \{\mathbf{n}\in\boldsymbol{\Lambda}^*\}}} \left(\liminf_{\substack{\mathbf{n}\to\infty \\ \{\mathbf{n}\in\boldsymbol{\Lambda}^*\}}}\right) \frac{T_{\mathbf{n},\mathbf{n}+\mathbf{n}^\alpha}}{\sqrt{2|\mathbf{n}|^\alpha\log|\mathbf{n}|}} = \sigma\sqrt{\frac{1-\alpha}{\beta}} \ \left(-\sigma\sqrt{\frac{1-\alpha}{\beta}}\right) \qquad a.s. \qquad (6.6)$$

*Conversely, if*

$$P\left(\limsup_{\substack{\mathbf{n}\to\infty \\ \{\mathbf{n}\in\boldsymbol{\Lambda}^*\}}} \frac{|T_{\mathbf{n},\mathbf{n}+\mathbf{n}^\alpha}|}{\sqrt{|\mathbf{n}|^\alpha\log|\mathbf{n}|}} < \infty\right) > 0, \qquad (6.7)$$

*then* (6.5) *holds,* $\mathrm{E}X = 0$ *and* (6.6) *holds with* $\sigma^2 = \operatorname{Var} X$.

The theorem tells us that the extreme limit points get closer and approach zero as $\beta$ increases, that is, the thinner the subsequence, the less wild are the observable oscillations.

**Proof of Theorem 6.1.** Since we consider subsequences, there are no gaps to fill. The only thing to check for the sufficiency is the independence for the lower bound, which is "immediate" since the subsequence here is sparser than that in (5.14). Indeed,

$$(i+1)^{\beta/(1-\alpha)} = i^{\beta/(1-\alpha)} \left(1 + \frac{1}{i}\right)^{\beta/(1-\alpha)} > i^{\beta/(1-\alpha)} \left(1 + \frac{\beta}{1-\alpha} \frac{1}{i}\right)$$

$$= i^{\beta/(1-\alpha)} + i^{(\beta-1+\alpha)/(1-\alpha)} \geq i^{\beta/(1-\alpha)} + i^{(\beta\alpha)/(1-\alpha)}.$$

The converse follows as in Theorem 2.1, so there is nothing more to prove there. $\qquad\square$

**Remark 6.1.** An even closer inspection of the proofs shows that, in fact, Theorem 6.1 remains true with $i^{\beta/(1-\alpha)}$ replaced by $(\frac{i}{\log i})^{\beta/(1-\alpha)}$ or even by $i^{\beta_1/(1-\alpha)}(\log i)^{\beta_2}$ for any $\beta_1 > 1$ and $\beta_2 \in \mathbb{R}$, and, more generally, by $i^{\beta/(1-\alpha)}\ell(i)$, where $\beta > 1$ and $\ell$ is a slowly varying function.

**Remark 6.2.** The results above show that the set of limit points of $\frac{T_{\mathbf{n},\mathbf{n}+\mathbf{n}^\alpha}}{\sqrt{2|\mathbf{n}|^\alpha \log |\mathbf{n}|}}$ is given by the whole interval $[-\sigma\sqrt{1-\alpha}, \sigma\sqrt{1-\alpha}]$.

**Remark 6.3.** Theorem 6.1 also holds for the case $d = 1$, thereby providing an extension to subsequences of Lai's original result.

As a special case, we mention the following result which tells us what happens if we consider the points along the diagonal. We leave the details to the reader.

**Theorem 6.2.** *Under the assumptions of Theorem 6.1,*

$$\limsup_{\substack{\mathbf{n}\to\infty \\ \{\mathbf{n}=(i^{\beta/(1-\alpha)},\dots,i^{\beta/(1-\alpha)})\}}} \left(\liminf_{\substack{\mathbf{n}\to\infty \\ \{\mathbf{n}=(i^{\beta/(1-\alpha)},\dots,i^{\beta/(1-\alpha)})\}}}\right) \frac{T_{\mathbf{n},\mathbf{n}+\mathbf{n}^\alpha}}{\sqrt{2|\mathbf{n}|^\alpha \log |\mathbf{n}|}}$$

$$= \sigma\sqrt{\frac{1-\alpha}{d\beta}} \quad \left(-\sigma\sqrt{\frac{1-\alpha}{d\beta}}\right) \qquad a.s.$$

*The converse is as before.*

# 7. Maximal windows

An LSL for the sequence of *maximal* windows or delayed sums is now easily attainable with the aid of the Lévy inequalities; see Lai [15] for the case $d = 1$.

**Theorem 7.1.** *Suppose that* $\{X_{\mathbf{k}}, \mathbf{k} \in \mathbf{Z}_+^d\}$ *are i.i.d. random variables with mean 0 and finite variance* $\sigma^2$, *and set* $S_{\mathbf{n}} = \sum_{\mathbf{k} \leq \mathbf{n}} X_{\mathbf{k}}$, $\mathbf{n} \in \mathbf{Z}_+^d$. *If* (2.1) *holds, then*

$$\limsup_{\mathbf{n} \to \infty} \frac{\max_{\mathbf{0} \leq \mathbf{k} \leq \mathbf{n}^{\alpha}} T_{\mathbf{n}, \mathbf{n} + \mathbf{k}}}{\sqrt{2 |\mathbf{n}|^{\alpha} \log |\mathbf{n}|}} = \sigma \sqrt{1 - \alpha} \qquad a.s. \tag{7.1}$$

*Conversely, if*

$$P\left( \limsup_{\mathbf{n} \to \infty} \frac{\max_{\mathbf{0} \leq \mathbf{k} \leq \mathbf{n}^{\alpha}} |T_{\mathbf{n}, \mathbf{n} + \mathbf{k}}|}{\sqrt{|\mathbf{n}|^{\alpha} \log |\mathbf{n}|}} < \infty \right) > 0, \tag{7.2}$$

*then* (2.1) *holds,* $\mathbf{E} X = 0$ *and* (7.1) *holds with* $\sigma^2 = \operatorname{Var} X$.

**Proof.** Since

$$|T_{\mathbf{n}, \mathbf{n} + \mathbf{n}^{\alpha}}| \leq \max_{\mathbf{0} \leq \mathbf{k} \leq \mathbf{n}^{\alpha}} |T_{\mathbf{n}, \mathbf{n} + \mathbf{k}}|,$$

the only implication that requires a proof is (2.1) $\Longrightarrow$ (7.1) and this is achieved with the aid of the Lévy inequality Lemma 3.1(b), according to which

$$P\left( \frac{\max_{\mathbf{0} \leq \mathbf{k} \leq \mathbf{n}^{\alpha}} |T_{\mathbf{n}, \mathbf{n} + \mathbf{k}}|}{\sqrt{|\mathbf{n}|^{\alpha} \log |\mathbf{n}|}} > \varepsilon \right) = P\left( \max_{\mathbf{0} \leq \mathbf{k} \leq \mathbf{n}^{\alpha}} |T_{\mathbf{n}, \mathbf{n} + \mathbf{k}}| > \varepsilon \sqrt{|\mathbf{n}|^{\alpha} \log |\mathbf{n}|} \right)$$

$$\leq 2^d P(|T_{\mathbf{n}, \mathbf{n} + \mathbf{n}^{\alpha}}| > \varepsilon \sqrt{|\mathbf{n}|^{\alpha} \log |\mathbf{n}|} - d \sqrt{2 \operatorname{Var}(T_{\mathbf{n}, \mathbf{n} + \mathbf{n}^{\alpha}})})$$

$$\leq 2^d P(|T_{\mathbf{n}, \mathbf{n} + \mathbf{n}^{\alpha}}| > \varepsilon (1 - \delta) \sqrt{|\mathbf{n}|^{\alpha} \log |\mathbf{n}|})$$

for any $\delta > 0$, provided $\mathbf{n}$ is sufficiently large. $\qquad \square$

The LSL for subsequences carries over similarly.

**Theorem 7.2.** *Suppose that* $\{X_{\mathbf{k}}, \mathbf{k} \in \mathbf{Z}_+^d\}$ *are i.i.d. random variables with mean 0 and finite variance* $\sigma^2$, *and set* $S_{\mathbf{n}} = \sum_{\mathbf{k} \leq \mathbf{n}} X_{\mathbf{k}}$, $\mathbf{n} \in \mathbf{Z}_+^d$. *If* (2.1) *holds, then, for* $\beta > 1$,

$$\limsup_{\substack{\mathbf{n} \to \infty \\ \{\mathbf{n}: |\mathbf{n}| = i^{\beta/(1-\alpha)}\}}} \frac{\max_{\substack{\mathbf{0} \leq \mathbf{k} \leq \mathbf{n}^{\alpha} \\ \{\mathbf{n}: |\mathbf{n}| = i^{\beta/(1-\alpha)}\}}} T_{\mathbf{n}, \mathbf{n} + \mathbf{k}}}{\sqrt{2 |\mathbf{n}|^{\alpha} \log |\mathbf{n}|}} = \sigma \sqrt{\frac{1 - \alpha}{\beta}} \qquad a.s. \tag{7.3}$$

*Conversely, if*

$$P\left(\limsup_{\substack{\mathbf{n}\to\infty\\ \{\mathbf{n}:|\mathbf{n}|=i^{\beta/(1-\alpha)}\}}} \frac{\max_{\substack{\mathbf{0}\le\mathbf{k}\le\mathbf{n}^\alpha\\ \{\mathbf{n}:|\mathbf{n}|=i^{\beta/(1-\alpha)}\}}} |T_{\mathbf{n},\mathbf{n}+\mathbf{k}}|}{\sqrt{|\mathbf{n}|^\alpha \log|\mathbf{n}|}} < \infty\right) > 0, \qquad (7.4)$$

*then* (2.1) *holds,* $\mathrm{E}X = 0$ *and* (7.3) *holds with* $\sigma^2 = \operatorname{Var} X$.

# 8. Additional results and remarks

In this closing section, we begin with some additional comments and remarks, after which we close with an LIL-type result derived via the delta method, which, in turn, is applied to the LIL as well as the LSL.

## 8.1. The set of limit points is independent of $d$

From the central limit theorem, we know that

$$\frac{T_{\mathbf{n},\mathbf{n}+\mathbf{n}^\alpha}}{\sqrt{|\mathbf{n}|^\alpha}} \xrightarrow{d} N(0,\sigma^2) \qquad \text{as } \mathbf{n}\to\infty.$$

Moreover, the windows are disjoint (at least asymptotically). This means that, heuristically,

$$\frac{T_{\mathbf{n},\mathbf{n}+\mathbf{n}^\alpha}}{\sqrt{|\mathbf{n}|^\alpha \log|\mathbf{n}|}} = \frac{T_{\mathbf{n},\mathbf{n}+\mathbf{n}^\alpha}}{\sqrt{|\mathbf{n}|^\alpha}} \cdot \frac{1}{\sqrt{\log|\mathbf{n}|}} \approx \frac{V_{\mathbf{n}}}{\sqrt{\log|\mathbf{n}|}},$$

where $\{V_{\mathbf{n}}, \mathbf{n}\in\mathbf{Z}_+^d\}$ are i.i.d. $N(0,\sigma^2)$-distributed random variables. This implies, among other things, that one might expect $\{\max_{\mathbf{n}} \frac{T_{\mathbf{n},\mathbf{n}+\mathbf{n}^\alpha}}{\sqrt{|\mathbf{n}|^\alpha \log|\mathbf{n}|}}, \mathbf{n}\in\mathbf{Z}_+^d\}$ to share the asymptotics of $\{\max_{\mathbf{n}} \frac{V_{\mathbf{n}}}{\sqrt{\log|\mathbf{n}|}}, \mathbf{n}\in\mathbf{Z}_+^d\}$, which would mean that, asymptotically, we are dealing with a sequence of maxima of i.i.d. Gaussian random variables normalized by the logarithm of the *number* of them, which, of course, does not depend on the structure of the index set.

For an interesting reference in the present context – maxima of Gaussian random fields – we refer to [16].

## 8.2. On the choice of subsequences

In the course of the proofs, we have seen that the double- and triple-primed sums are negligible and that the size of the oscillations depends on the primed sums along a suitably selected subsequence. In the proofs of the LIL for sums, the exponential bounds are exponentials of *iterated* logarithms, that is, powers of logarithms. In order to obtain a convergent sum (for the upper bound), the natural choices of subsequences are geometrically increasing ones.

In our context, the exponential bounds are exponentials of *single* logarithms, that is, powers. In order to obtain a convergent sum (for the upper bound), the natural choices of subsequences are polynomially increasing ones.

We also recall that in connection with Theorem 6.1, we observed that the oscillations become less wild as the subsequences get thinner.

## 8.3. Functional LSL's

A further project might be to consider possible functional or Strassen versions of our results. For references in the one-dimensional case, see, for example, [2, 6].

## 8.4. The LIL, the LSL and the delta method

The following application of the so-called delta method to the law of the iterated logarithm might not be new. However, we have never seen it in the literature (cf., e.g., [11], page 349 in connection with the central limit theorem).

**Theorem 8.1.** *Suppose that $\{U_n, n \geq 1\}$ is a sequence of random variables, for which there exist positive sequences $\{a_n, n \geq 1\}$ and $\{b_n, n \geq 1\}$ tending to infinity as $n \to \infty$, such that*

$$\frac{U_n}{b_n} \overset{\text{a.s.}}{\to} \mu \qquad \text{as } n \to \infty \ (0 \leq \mu < \infty)$$

*and*

$$\limsup_{n \to \infty} \left( \liminf_{n \to \infty} \right) a_n \left( \frac{U_n}{b_n} - \mu \right) = 1 \ (-1) \qquad a.s. \tag{8.1}$$

(i) *If $g$ is continuously differentiable at (in a neighbourhood of) $\mu$ and $g'(\mu) \neq 0$, then*

$$\limsup_{n \to \infty} \left( \liminf_{n \to \infty} \right) a_n \left( g\left( \frac{U_n}{b_n} \right) - g(\mu) \right) = +|g'(\mu)| \ (-|g'(\mu)|).$$

(ii) *If $g$ is $m \geq 2$ times continuously differentiable at (in a neighbourhood of) $\mu$, $g^{(k)}(\mu) = 0$, $k = 1, 2, \ldots, m-1$, and $g^{(m)}(\mu) > 0$, say, then*

$$\limsup_{n \to \infty} \left( \liminf_{n \to \infty} \right) a_n^{m/2} \left( g\left( \frac{U_n}{b_n} \right) - g(\mu) \right) = +\frac{1}{m!} g^{(m)}(\mu) \ \left( -\frac{1}{m!} g^{(m)}(\mu) \right),$$

*if $m$ is odd. If, on the other hand, $m$ is even, then*

$$\limsup_{n \to \infty} a_n^{m/2} \left( g\left( \frac{U_n}{b_n} \right) - g(\mu) \right) = +\frac{1}{m!} g^{(m)}(\mu).$$

**Proof.** (i) By Taylor expansion,

$$g\left(\frac{U_n}{b_n}\right) = g(\mu) + \left(\frac{U_n}{b_n} - \mu\right)g'(\theta_n),$$

where $|\theta_n - \mu| \leq |\frac{U_n}{b_n} - \mu|$, so that

$$a_n\left(g\left(\frac{U_n}{b_n}\right) - g(\mu)\right) = a_n\left(\frac{U_n}{b_n} - \mu\right)g'(\theta_n),$$

and (i) follows.

(ii) Using second-order Taylor expansion, and recalling that $g'(\mu) = 0$, we similarly obtain

$$a_n^2\left(g\left(\frac{U_n}{b_n}\right) - g(\mu)\right) = \frac{1}{2}\left(a_n\left(\frac{U_n}{b_n} - \mu\right)\right)^2 g''(\theta_n),$$

which establishes (ii) in the case $m = 2$ and similarly in the other cases with higher order Taylor expansions. Note that one loses the sign of $\left(\frac{U_n}{b_n} - \mu\right)$ when $m$ is even. $\qquad\square$

***Remark 8.1.*** In case the limit points in (8.1) are dense in $[-1, 1]$, the liminf in case (ii) for $g^{(m)}(\mu) > 0$ with $m$ even is equal to zero.

As immediate corollaries, we obtain the following results related to the Hartman–Wintner LIL ([13] or, e.g., [11], Theorem 8.1.2) and Lai's LSL for delayed sums [15], respectively.

***Corollary 8.1.*** *Suppose that $X, X_1, X_2, \ldots$ are i.i.d. random variables with mean $\mu$ and finite variance $\sigma^2$, and set $S_n = \sum_{k=1}^n X_k$, $n \geq 1$.*

(i) *If $g$ is continuously differentiable at (in a neighbourhood of) $\mu$ and $g'(\mu) \neq 0$, then*

$$\limsup_{n\to\infty}\left(\liminf_{n\to\infty}\right)\sqrt{\frac{n}{\log\log n}}\left(g\left(\frac{S_n}{n}\right) - g(\mu)\right) = +\sigma\sqrt{2}|g'(\mu)| \ (-\sigma\sqrt{2}|g'(\mu)|).$$

(ii) *If $g$ is $m \geq 2$ times continuously differentiable at (in a neighbourhood of) $\mu$, $g^{(k)}(\mu) = 0$, $k = 1, 2, \ldots, m-1$, and $g^{(m)}(\mu) > 0$, say, then*

$$\limsup_{n\to\infty}\left(\liminf_{n\to\infty}\right)\left(\frac{n}{\log\log n}\right)^{m/2}\left(g\left(\frac{S_n}{n}\right) - g(\mu)\right)$$

$$= +\frac{1}{m!}(\sigma\sqrt{2})^m g^{(m)}(\mu) \ \left(-\frac{1}{m!}(\sigma\sqrt{2})^m g^{(m)}(\mu)\right),$$

*if $m$ is odd. If, on the other hand, $m$ is even, then*

$$\limsup_{n\to\infty}\left(\frac{n}{\log\log n}\right)^{m/2}\left(g\left(\frac{S_n}{n}\right) - g(\mu)\right) = +\frac{1}{m!}(\sigma\sqrt{2})^m g^{(m)}(\mu)$$

*and*

$$\liminf_{n\to\infty}\left(\frac{n}{\log\log n}\right)^{m/2}\left(g\left(\frac{S_n}{n}\right)-g(\mu)\right)=0.$$

**Proof.** Use Theorem 8.1 with $U_n=S_n$, $b_n=n$, $\mu=\mathrm{E}X$ and $a_n=\sqrt{\frac{n}{2\sigma^2\log\log n}}$, and apply the strong law of large numbers and the Hartman–Wintner LIL in order to verify that the assumptions are fulfilled. An appeal to Remark 8.1 completes the proof. $\qquad\square$

**Corollary 8.2.** *Let $X, X_1, X_2, \ldots$ be i.i.d. random variables with mean 0 and variance 1 and suppose that $\mathrm{E}|X|^{2/\alpha}(\log^+|X|)^{-1/\alpha}<\infty$. Finally, set $T_{n,n+n^\alpha}=\sum_{k=n+1}^{n+n^\alpha}X_k$.*

(i) *If $g$ is continuously differentiable at (in a neighbourhood of) 0 and $g'(0)\neq0$, then*

$$\limsup_{n\to\infty}\left(\liminf_{n\to\infty}\right)\sqrt{\frac{n^\alpha}{\log n}}\left(g\left(\frac{T_{n,n+n^\alpha}}{n^\alpha}\right)-g(0)\right)=\sqrt{2}|g'(0)|\ (-\sqrt{2}|g'(0)|).$$

(ii) *If $g$ is $m\geq2$ times continuously differentiable at (in a neighbourhood of) $\mu$, $g^{(k)}(\mu)=0$, $k=1,2,\ldots,m-1$, and $g^{(m)}(\mu)>0$, say, then*

$$\limsup_{n\to\infty}\left(\liminf_{n\to\infty}\right)\left(\frac{n^\alpha}{\log n}\right)^{m/2}\left(g\left(\frac{T_{n,n+n^\alpha}}{n^\alpha}\right)-g(0)\right)$$

$$=+\frac{1}{m!}2^{m/2}g^{(m)}(0)\quad\left(-\frac{1}{m!}2^{m/2}g^{(m)}(0)\right),$$

*if $m$ is odd. If, on the other hand, $m$ is even, then*

$$\limsup_{n\to\infty}\left(\frac{n^\alpha}{\log n}\right)^{m/2}\left(g\left(\frac{T_{n,n+n^\alpha}}{n^\alpha}\right)-g(0)\right)=+\frac{1}{m!}2^{m/2}g^{(m)}(0)$$

*and*

$$\liminf_{n\to\infty}\left(\frac{n^\alpha}{\log n}\right)^{m/2}\left(g\left(\frac{T_{n,n+n^\alpha}}{n^\alpha}\right)-g(0)\right)=0.$$

**Proof.** Use Theorem 8.1 with $U_n=T_{n,n+n^\alpha}$, $b_n=n^\alpha$, $\mu=0$ and $a_n=\sqrt{\frac{n^\alpha}{2\log n}}$, and apply Chow's strong law [5] and Lai's LSL [15] in order to check the assumptions. For the case where $m$ is even, use results of this paper together with Remark 8.1. $\qquad\square$

# Acknowledgements

The basis of this paper originates with the first author's visit to the Department of Number and Probability Theory at the University of Ulm. Allan Gut wishes to thank his coauthor, Professor Stadtmüller, for providing a wonderful stay, for two most stimulating weeks and for his considerate and generous hospitality. He also wishes to thank the University of Ulm for financial support.